# Dimensionality reduction based on Distance Preservation to Local Mean (DPLM) for SPD matrices and its application in BCI

Alireza Davoudi[a], Saeed Shiry Ghidary[a], Khadijeh Sadatnejad[a]

[a] Computer engineering and information technology department, Amirkabir university of technology, Tehran, Iran

**Abstract**

In this paper, we propose a nonlinear dimensionality reduction algorithm for the manifold of Symmetric Positive Definite (SPD) matrices that considers the geometry of SPD matrices and provides a low dimensional representation of the manifold with high class discrimination. The proposed algorithm, tries to preserve the local structure of the data by preserving distance to local mean (DPLM) and also provides an implicit projection matrix. DPLM is linear in terms of the number of training samples and may use the label information when they're available in order to performance improvement in classification tasks. We performed several experiments on the multi-class dataset IIa from BCI competition IV. The results show that our approach as dimensionality reduction technique - leads to superior results in comparison with other competitor in the related literature because of its robustness against outliers. The experiments confirm that the combination of DPLM with FGMDM as the classifier leads to the state of the art performance on this dataset.

## 1. Introduction

Motor Imagery (MI) based Brain Computer Interface (BCI) relies on using electrical activity of brain - which is usually Electroencephalography (EEG) signals - to decode human motor intentions. This kind of BCI is very useful in rehabilitating sensorimotor functions in patients who have severe motor disabilities [1-5].

The most common approach for feature extraction from EEG signals in MI paradigm is Common Spatial Patterns (CSP) [6]. CSP projects multichannel EEG signals into a more discriminative subspace, which can have a lower dimension. Therefore, we can consider CSP as a nonlinear dimensionality reduction technique [7]. CSP is originally designed for two class problems but there are multiclass extensions available [8, 9]. The performance of CSP is highly dependent to the number of electrodes and also their locations [10].

Using covariance matrices for describing data has attracted more attention recently in several fields specially, computer vision [11] and brain computer interfacing [12]. Providing a compact and informative representation are the motivations for using covariance matrices as descriptor. However, in many applications that rely on using covariance matrices as descriptors, such as MI BCI and computer aided diagnosis systems, insufficient training samples might cause the curse of dimensionality problem. Therefore, it is critical to use a dimensionality reduction technique.

Covariance matrices lie in the space of Symmetric Positive Definite (SPD) matrices, which can be formulated as a Riemannian manifold. There are three major approaches for reducing the dimensionality of SPD matrices. The first approach is to map the SPD matrices into tangent spaces and then use Euclidean methods for dimensionality reduction in resulting Euclidean space. Principal Geodesic Analysis (PGA) [13] is based on this approach and provides a generalization of Principal Component Analysis (PCA) to Riemannian manifolds. It tries to find a tangent space that maximizes the variability of mapped data points. However, PGA is equivalent to mapping the data to the tangent space of the geometric mean of the data [14]. Goh et.al. [15], provide the generalization of three local Nonlinear Dimensionality Reduction (NLDR) algorithms including Locally Linear Embedding (LLE), Laplacian Eigenmaps (LE) and Hessian LLE (HLLE) to Riemannian manifold. They provide an embedding to lower dimensional space based on Riemannian geometry, for example using Karcher mean and Riemannian logarithmic map. Since, these NLDR algorithms doesn't provide an implicit mapping, they can't be used in non-transductive scenarios.

The second approach for dimensionality reduction over Riemannian manifolds relies on kernel approach, which tries to embed SPD matrices in a Reproducing Kernel Hilbert Space (RKHS) and then perform dimensionality reduction using existing kernel based methods such as KPCA. Goh et. Al. [16] proposed a kernel function that embeds the SPD manifold into a Euclidean space. This mapping can be considered as mapping the SPD manifold into the tangent space at identity matrix [14]. Barachant et. al. [21] proposed a kernel that is very similar to the one proposed by Goh et. al. [16] with the difference at base point of tangent space. They choose the Karcher mean as the base point of the tangent space. Jayasumana et. al. [17] proposed a positive definite Gaussian RBF kernel, which embeds the SPD manifold into an RKHS. In [18] a Riemannian pseudo kernel is proposed, which is positive definite under certain conditions. In [7] a kernel is proposed by manipulating the indefinite isometric kernel using a geometry preserving conformal transform. As the dimension of SPD matrices grows, the computational cost of kernel approaches highly increases.

The third approach for overcoming the problem of high dimensionality in SPD manifolds is mapping from a high dimensional SPD manifold to a lower dimensional one while the geometry of SPD manifolds is preserved. This kind of dimensionality reduction has two important properties: 1) it directly works on the original manifold to learn a mapping [14] and 2) the resulting low dimensional manifold can be used as the input to existing SPD based algorithms. The only work of this kind is Harandi's work [14]. They learn a mapping that maximize the geodesic distances between inter-class samples and simultaneously minimize the distances between intra-class samples and it is done via an optimization on Grassmann manifolds. This method has two drawbacks: 1) As all of the pairs of samples have the same weight in the optimization formula, the result is prone to outlier samples and 2) One cycle of the optimization process has an order of $O(N^2)$ where $N$ is the number of samples and thus it becomes intractable when the number of samples grows.

In this paper, we propose a dimensionality reduction algorithm for the space of SPD matrices, which can be considered as a method that belongs to the third approach. This algorithm tries to preserve the local structure of the data by distance preservation to local mean (DPLM), considers the geometry of SPD matrices, provides an implicit mapping and applies the supervised information for embedding to lower dimensional space. As we just use the neighbors of every sample in the optimization process, there is a good chance that outliers don't contribute in the optimization. One cycle of the optimization process has an order of $O(KN)$ where $N$ is the number of samples and $K$ is the number of the neighbors of every sample. Thus it is linear in terms of number of samples.

The rest of this paper is organized as follows: in section 2, we describe the mathematical preliminaries. Details of the proposed algorithm are presented in section 3. Section 4 reports our experiments on a BCI dataset and also a comparison with Harandi's method. We conclude the results in section 5.

## 2. Geometry of SPD matrices

A $d$-dimensional *topological manifold* is a connected paracompact Hausdorff space that is locally homeomorphic to the $d$-dimensional Euclidean space $\mathbb{R}^d$. A *differentiable manifold* is a topological manifold that has a globally defined differential structure. The *tangent space* at a point $p$ on a differentiable manifold is the plane tangent to the surface of the manifold at that point. A *Riemannian manifold* is a differentiable manifold with a *Riemannian metric*. The Riemannian metric of a manifold consists of all inner products on all of the tangent spaces of the manifold. We can measure the angle between two curves and also the length of a curve by using a Riemannian metric. The *geodesic distance* between two points on a Riemannian manifold is the length of the shortest path between these two points.

A real symmetric matrix $X \in \mathbb{R}^{d \times d}$ is a SPD matrix if and only if $v^T X v > 0$ for every non-zero $v \in \mathbb{R}^d$. The set of $d \times d$ SPD matrices is denoted by $Sym_d^+$, which is a differentiable manifold with a natural Riemannian structure. The most common Riemannian metric proposed on $Sym_d^+$ is Affine Invariant Riemannian Metric (AIRM)[19]:

$$\delta_g^2(X, Y) = \|log(XY^{-1})\|_{Frob}^2 \tag{1}$$

where $X, Y \in Sym_d^+$ and $\log(.)$ is matrix logarithm. AIRM is a true geodesic distance but is computationally expensive, especially when the dimension is larger than 20 [20]. Jensen-Bregman Log-det Divergence (JBLD) [20] has been proposed as a similarity measure for SPD matrices:

$$J(X, Y) = logdet\left(\frac{X+Y}{2}\right) - \frac{1}{2} logdet(XY) \tag{2}$$

where $X, Y \in Sym_d^+$. It has been shown that:

$$\delta_{ld}^2(X, Y) = \sqrt{J(X, Y)} \tag{3}$$

is a metric on $Sym_d^+$ [21] and has a lower computational cost compared to AIRM [20].

The mean of $N$ SPD matrices, which is also referred to as the geometric mean, is given by:

$$\varphi(P_1, \dots, P_N) = \underset{P \in Sym_d^+}{arg\,min} \sum_{i=1}^{N} \delta_g^2(P, P_i) \qquad (4)$$

This mean exist and is unique [22] but has no closed form expression. An iterative method for computing the geometric mean is proposed in [23].

## 3. The proposed algorithm

Given an SPD manifold M, we wish to find a lower dimensional representation of it at the space of SPD matrices by preserving the local structure, which is done by preserving distance to Local Mean (DPLM). To this end, we calculate the Riemannian mean of the K nearest neighbor of each training sample and try to find a projection matrix that preserves the distances between each of the K nearest neighbors and their means. This is illustrated in Fig. 1. As we just use the neighbors of each sample in the optimization process, there is a good chance that the outliers doesn't contribute in the optimization and therefore leads to robustness against outliers.

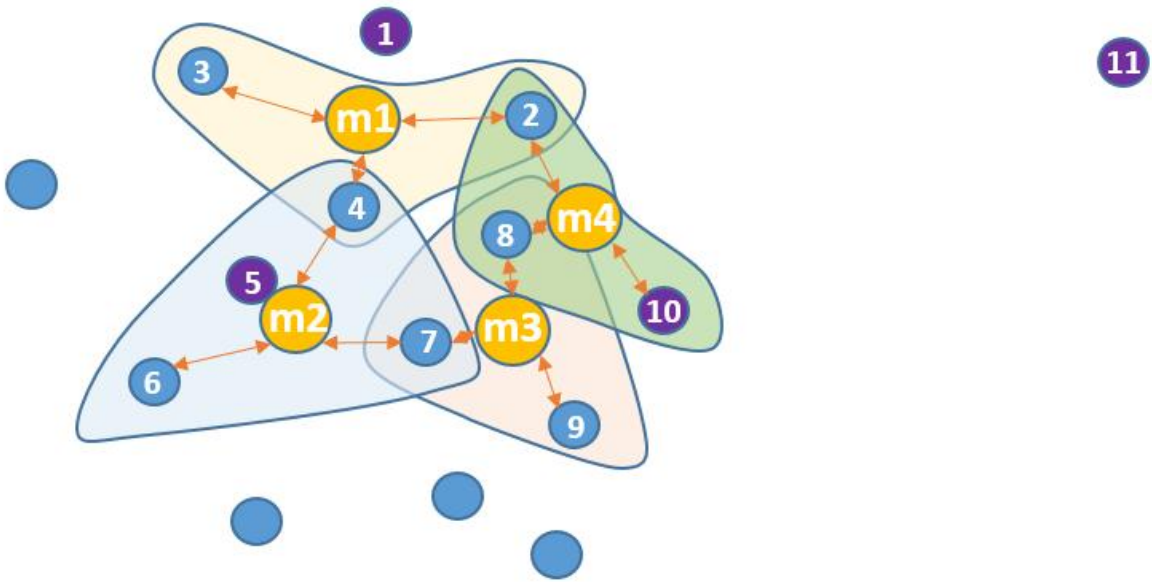

**Fig. 1.** An illustration of DPLM. m1 is the mean of 2,3 and 4, which are the three NNs of 1. m2 is the mean of 4,6 and 7, which are the three NNs of 5. m3 is the mean of 7,8 and 9, which are the three NNs of 10. m4 is the mean of 2,8 and 10, which are the three NNs of 11. The two sided arrows mean that we want to preserve the distance between the two samples. Note that sample 11 is an example of an outlier that doesn't contribute in the optimization, as it's not in the 3 NNs set of the other samples.

We now describe DPLM in more details. Suppose we have $N$ training samples $\{(X_1, y_1), (X_2, y_2), \ldots, (X_N, y_N)\}$ where $X_i \in Sym_n^+$ and $y_i$ is its corresponding class label. We aim to find a projection matrix $U$ that maps the samples to $Sym_m^+$ where $m < n$.

Suppose $N_i = \{X_{i,1}, X_{i,2}, \ldots, X_{i,k}\}, 1 \leq i \leq N$ is the set of $K$ nearest neighbors of sample $X_i$. For classification purpose, it's better to select the $K$ nearest neighbors that have the same label as that of $X_i$, i.e. $y_i$. The Riemannian mean of each set $N_i$ denoted by $\overline{N}_i$ is calculated using Eq. (4).

Next we find the projection matrix $U$ by solving the following optimization problem:

$$\min_{U \in \mathbb{R}^{n \times m}} H(U), \quad s.t \quad U^T U = I_m \tag{5}$$

where

$$H(U) = \sum_{i=1}^{N} \sum_{j=1}^{K} \left| \delta_{ld}^2(X_{i,j}, \overline{N}_i) - \delta_{ld}^2(U^T X_{i,j} U, U^T \overline{N}_i U) \right| \tag{6}$$

$U = [U_1, U_2, \ldots, U_m] \in \mathbb{R}^{n \times m}$ is an orthonormal matrix and $I_m$ is an $m \times m$ identity matrix. The constraint $U^T U = I_m$ is required to ensure that $U^T X U$ is positive definite for every $X \in Sym_n^+$ and therefore is a valid SPD matrix. This constraint also helps to avoid degeneracy in the optimization process.

The JBLD metric, which provides the same result as $\delta_g^2$ [21] up to a scale of $2\sqrt{2}$ [14], is an appropriate choice for Eq. (6) because of its lower computations.

Wen et. al. [24] proposed a method for optimization of problems with orthogonality constraints similar to Eq. (5). This constraint is non-convex and satisfying it during iterations is numerically expensive. Given a feasible point $X$, and its corresponding gradient $G$, the new trial point with the step size $\tau \geq 0$ is obtained by the Crank-Nicolson-like scheme:

$$Y(\tau) = \left(I + \frac{\tau}{2} A\right)^{-1} \left(I + \frac{\tau}{2} A\right) X \tag{7}$$

where $A := G X^T - X G^T$ is a skew symmetric matrix. Because of some nice properties of $Y(\tau)$, the $\{Y(\tau)\}_{\tau \geq 0}$ is a descent path. Therefore, to guarantee the convergence, a curvilinear search is applied to find a proper step size ($\tau$).

To use Wen's approach for solving Eq. (5) the gradient of Eq. (6) with respect to $U$ is computed as follows:

$$\frac{\partial H(U)}{\partial U} = -\sum_{i=1}^{N}\sum_{j=1}^{K}\left[sgn\left(\delta_{ld}^2(X_{i,j},\bar{N}_i) - \delta_{ld}^2(U^TX_{i,j}U, U^T\bar{N}_iU)\right)\right.$$
$$\times\left((X_{i,j}+\bar{N}_i)U\left(U^T\frac{X_{i,j}+\bar{N}_i}{2}U\right)^{-1} - X_{i,j}U(U^TX_{i,j}U)^{-1}\right. \quad (8)$$
$$\left.\left.- \bar{N}_iU(U^T\bar{N}_iU)^{-1}\right)\right]$$

with the following prior knowledge [14]:

$$\frac{\partial\delta_{ld}^2(U^TXU, U^TYU)}{\partial U}$$
$$= (X+Y)U\left(U^T\frac{X+Y}{2}U\right)^{-1} - XU(U^TXU)^{-1} - YU(U^TYU)^{-1} \quad (9)$$

where $sgn(.)$ is the sign function.

The lower dimensional representation of a new sample $X_{new}$ in $Sym_m^+$ space can be computed as:

$$X'_{new} = U^TX_{new}U \in Sym_m^+ \quad (10)$$

where $X'_{new}$ is the representation of $X_{new}$ in $Sym_m^+$ space.

## 4. Evaluations

To evaluate our dimensionality reduction method, we use the Dataset IIa of BCI competition IV [25]. We use the proposed method in conjunction with the two simple but powerful classifiers: Minimum Distance to Mean (MDM) [12] and Filter Geodesic MDM (FGMDM) [12]. In MDM, first the Riemannian mean of each class is calculated, then a test sample is assigned to the class that has the shortest distance to its mean. Since MDM is not robust to noise, it is suggested to perform some filtering over the data before applying MDM [26]. FGMDM first tries to find a set of filters by applying an extension of Fisher Linear Discriminant Analysis (FLDA) named Fisher Geodesic Discriminant Analysis (FGDA) [12] and then apply these filters to data using a geodesic filtering approach that result in a set of SPD matrices with the same dimensionality as initial data. At last MDM will be used for classification.

### 4.1. Dataset description and preprocessing

Dataset IIa of BCI competition IV consists of EEG signals of 9 subjects. Each subject performs four kinds of motor imagery tasks (foot, tongue, right hand and left hand). Twenty-two electrodes are used for recording EEG signals. For each subject two sessions on two different days are recorded. One of them is for training and the other is for testing. Each session consists of 288 trials, i.e. 72 trials per class. This is a classification problem where we need to assign each testing sample to one of the four classes.

Before performing feature extraction, some preprocessing steps like the band-pass filtering and time interval selection are necessary. It has been shown that these pre-processing steps have a large impact on the final classification performance [26]. Here we propose a simple yet effective band-pass filter and time window selection algorithm that is based on cross-validation.

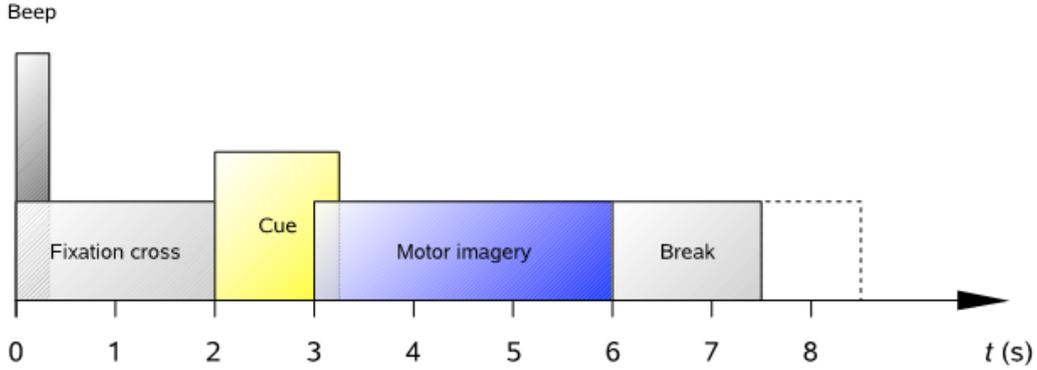

**Fig. 2.** Time schema of one trial of dataset IIa from BCI competition IV [25]

The timing schema of one trial is illustrated in Fig. 2. At the beginning of a trial a fixation cross appears on the screen and lasts for 2 seconds. At t = 2s a cue corresponding to one of the four classes appears and stays on screen for 1.25 seconds. After that the user starts the motor imagery task until the fixation cross disappears from the screen.

To find the efficient band-pass filter and time window for every subject, we use their training data. Next, we consider all of the possible windows with lengths from 1s to 4s with step 0.25s (i.e. 1, 1.25, 1.5, 1.75,…) from the starting point t = 3s, together with the set of four band-pass filters with the lower bands 5Hz and 8Hz and the higher bands 30Hz and 35Hz. Therefore, we have 180 test cases for evaluation. Note that each test case consists of a time window and a band-pass filter, i.e. we search for the efficient combination of time window and band pass filter. The evaluation is performed using 10-fold cross-validation and MDM. At last the mean of the K test cases with highest cross-validation accuracies is selected as the efficient time window and band pass filter for the corresponding subject. Here we selected K = 10.

Suppose $X \in \mathbb{R}^{22 \times N}$ is the band-passed EEG signal extracted from a desired time window of a trial. The covariance matrix of $X$ is calculated as follows:

$$C = \frac{1}{N-1} XX^T \qquad (11)$$

In this paper, we use $C$ as a descriptor for $X$.

Table 1 illustrates the results of using the proposed approach for time window and band-pass filter selection compared to when a fixed time window and band-pass filter is selected for all the subjects. As you can see, both MDM and FGMDM leads to higher kappa values where the proposed approach for preprocessing is used. Note that we selected the fixed band-pass filter as 8-35 Hz and the fixed time window as 3.75s to 5.75s. These are common values, which have been used in many studies [7, 27-29].

Table 1. Result of MDM and FGMDM after applying proposed preprocessing for selecting time window and band-pass filter. We compare it with the results achieved over fixed time window and band-pass filter (showed by a "fixed" postfix) in terms of kappa value.

|    | Time Window (seconds) | | Band-pass filter (Hz) | | MDM | FGMDM | MDM_fixed | FGMDM_fixed |
|----|-------|-------|-------|-------|------|------|------|------|
|    | start | end   | Lower | upper |      |      |      |      |
| S1 | 3.1   | 5.848 | 8     | 32    | **0.73** | 0.72 | 0.71 | 0.69 |
| S2 | 3.1   | 5.848 | 8     | 31    | 0.45 | **0.50** | 0.39 | 0.35 |
| S3 | 3     | 5.248 | 8     | 32    | 0.56 | **0.64** | 0.52 | 0.60 |
| S4 | 3     | 5.9   | 6     | 32    | **0.50** | 0.38 | 0.43 | 0.28 |
| S5 | 3     | 4.948 | 8     | 35    | 0.23 | **0.28** | 0.12 | 0.21 |
| S6 | 3.2   | 6     | 8     | 32    | 0.24 | **0.34** | 0.21 | 0.30 |
| S7 | 3     | 5.8   | 6     | 33    | 0.39 | **0.64** | 0.32 | 0.46 |
| S8 | 3.048 | 5.9   | 8     | 32    | 0.60 | **0.68** | 0.55 | 0.62 |
| S9 | 3     | 4.1   | 6     | 33    | 0.57 | **0.75** | 0.55 | 0.53 |
|    | **Average Kappa** | | | | 0.47 | **0.55** | 0.42 | 0.45 |

Hereafter, we use the covariance descriptors of the EEG data that have been preprocessed using the proposed method.

### 4.2. Experiments

To the best of our knowledge, the only competitor for embedding an SPD manifold into a lower dimensional symmetric positive definite space is Harandi's method [14]. Table 2 reports the results of MDM classifier for comparison between DPLM and Harandi's method, where 22 dimensional SPD matrices reduced to two up to 12 dimensional SPD matrices. The significant superiority of the proposed method in comparison with Harandi's method has been confirmed using the two-tailed Wilcoxson signed-rank statistical test, which results in p-value 0.01 (<0.05). To apply this test, we supposed that every pair of subject and dimension is a separated domain.

Table 2. A comparison between DPLM and Harandi's method in 2 up to 12 dimensional SPD space using the MDM classifier

|      | **DPLM** | | | | | | **Harandi's method** | | | | | |
|------|------|------|------|------|------|------|------|------|------|------|------|------|
| Dim. | 2    | 4    | 6    | 8    | 10   | 12   | 2    | 4    | 6    | 8    | 10   | 12   |
| S1   | 0.46 | 0.58 | 0.63 | 0.64 | 0.68 | 0.70 | 0.36 | 0.50 | 0.63 | 0.62 | 0.65 | 0.68 |
| S2   | 0.31 | 0.21 | 0.37 | 0.43 | 0.45 | 0.44 | 0.19 | 0.30 | 0.32 | 0.43 | 0.44 | 0.44 |
| S3   | 0.48 | 0.58 | 0.59 | 0.62 | 0.58 | 0.62 | 0.47 | 0.53 | 0.57 | 0.57 | 0.56 | 0.57 |
| S4   | 0.22 | 0.33 | 0.38 | 0.39 | 0.41 | 0.46 | 0.18 | 0.25 | 0.35 | 0.40 | 0.42 | 0.41 |
| S5   | 0.10 | 0.16 | 0.19 | 0.20 | 0.18 | 0.19 | 0.04 | 0.14 | 0.19 | 0.16 | 0.21 | 0.23 |

|     |      |      |      |      |      |      |      |      |      |      |      |
|-----|------|------|------|------|------|------|------|------|------|------|------|
| S6  | 0.15 | 0.17 | 0.19 | 0.20 | 0.18 | 0.19 | 0.15 | 0.15 | 0.14 | 0.18 | 0.14 | 0.19 |
| S7  | 0.19 | 0.30 | 0.35 | 0.27 | 0.32 | 0.36 | 0.21 | 0.18 | 0.24 | 0.29 | 0.24 | 0.30 |
| S8  | 0.49 | 0.59 | 0.62 | 0.64 | 0.63 | 0.61 | 0.34 | 0.58 | 0.58 | 0.57 | 0.61 | 0.56 |
| S9  | 0.45 | 0.47 | 0.52 | 0.53 | 0.55 | 0.55 | 0.43 | 0.44 | 0.49 | 0.51 | 0.53 | 0.56 |
| Avg.| **0.31** | **0.38** | **0.43** | **0.44** | **0.45** | **0.47** | **0.26** | **0.34** | **0.39** | **0.41** | **0.42** | **0.44** |

In Table 3, the performance of the proposed method, Harandi's method, FGMDM and the three first winners of BCI Competition IV on Dataset IIa are compared in terms of Kappa value. We used FGMDM for classifying the data in lower dimensional space. We should note that the proper dimension for each subject were achieved using 10-fold cross-validation. As it can be seen, the proposed approach when combined with FGMDM, achieves the highest mean Kappa, which is three percent better than the Harandi's mean Kappa.

**Table 3.** The performance of DPLM, Harandi's method, FGMDM and the three first winners of BCI Competition IV on Dataset IIa in terms of Kappa value. (The values in the parentheses show the dimensionality of the resulting space)

|          | S1       | S2       | S3       | S4       | S5       | S6       | S7       | S8       | S9       | Mean Kappa |
|----------|----------|----------|----------|----------|----------|----------|----------|----------|----------|------------|
| DPLM     | **0.74**(15) | 0.48 (18) | 0.72 (14) | 0.47 (18) | 0.28 (20) | 0.33 (20) | 0.63 (20) | 0.72 (17) | **0.76** (20) | 0.571 |
| 1st      | 0.68     | 0.42     | **0.75** | 0.48     | **0.40** | 0.27     | **0.77** | 0.75     | 0.61     | 0.57       |
| FGMDM    | 0.72     | **0.50** | 0.64     | 0.38     | 0.28     | **0.34** | 0.64     | 0.68     | 0.75     | 0.55       |
| Harandi  | 0.68 (19) | 0.43 (19) | 0.72 (19) | 0.46 (16) | 0.25 (15) | 0.31 (19) | 0.64 (20) | 0.63 (20) | 0.68 (20) | 0.54 |
| 2nd      | 0.69     | 0.34     | 0.71     | 0.44     | 0.16     | 0.21     | 0.66     | 0.73     | 0.69     | **0.52**   |
| 3rd      | 0.38     | 0.18     | 0.48     | 0.33     | 0.07     | 0.14     | 0.29     | 0.49     | 0.44     | **0.31**   |

The best performances for DPLM and Harandi's method in terms of Kappa value are reported in Table 4. As it can be seen, the mean Kappa of DPLM has superiority in comparison with Harandi's mean Kappa (three percent). To have a better comparison between these two methods, for each subject we use the lower dimension of the two methods in Table 4 and report the resulting performances in terms of Kappa value in Table 5. The results show that DPLM performs significantly better than Harandi's method.

**Table 4.** Comparison between the best performances for DPLM and Harandi's method in terms of Kappa value. (The values in the parentheses are dimensions)

|         | S1 | S2 | S3 | S4 | S5 | S6 | S7 | S8 | S9 | Mean Kappa |
|---------|----|----|----|----|----|----|----|----|----|------------|
| DPLM    | **0.75** (19) | **0.49** (14) | **0.76** (12) | 0.49 (14) | **0.34** (17) | **0.36** (15) | **0.68** (15) | **0.76** (14) | **0.76** (20) | 0.60 |
| Harandi | 0.73 (20) | 0.46 (20) | 0.75 (15) | 0.49 **(9)** | 0.25 (14) | 0.34 (14) | 0.65 (18) | 0.71 (18) | 0.74 (20) | 0.57 |

**Table 5.** Performance of DPLM and Harandi's method when for every subject the lower dimensions that are reported in Table 4 is used.

|        | S1   | S2   | S3   | S4   | S5   | S6   | S7   | S8   | S9   | Mean Kappa |
|--------|------|------|------|------|------|------|------|------|------|------------|
| **DPLM**   | **0.75** | **0.49** | **0.76** | 0.46 | **0.26** | 0.33 | **0.68** | **0.76** | **0.76** | **0.58** |
| **Harandi** | 0.68 | 0.38 | 0.72 | **0.49** | 0.25 | **0.34** | 0.61 | 0.67 | 0.74 | 0.54 |

In Fig. 3 (a), the execution time of DPLM for different values of K and the Harandi's algorithm versus the size of the sample set is illustrated. As we noted before, our method is linear in terms of the number of samples ($N$), however Harandi's method has an order of $O(N^2)$. Therefore, as it can be seen from Fig. 3 (a), DPLM performs much faster than Harandi's method when the number of samples grows. We have also compared the execution time of these two methods versus the number of dimensions that is illustrated in Fig. 3 (b) and as it can be seen, DPLM is much faster.

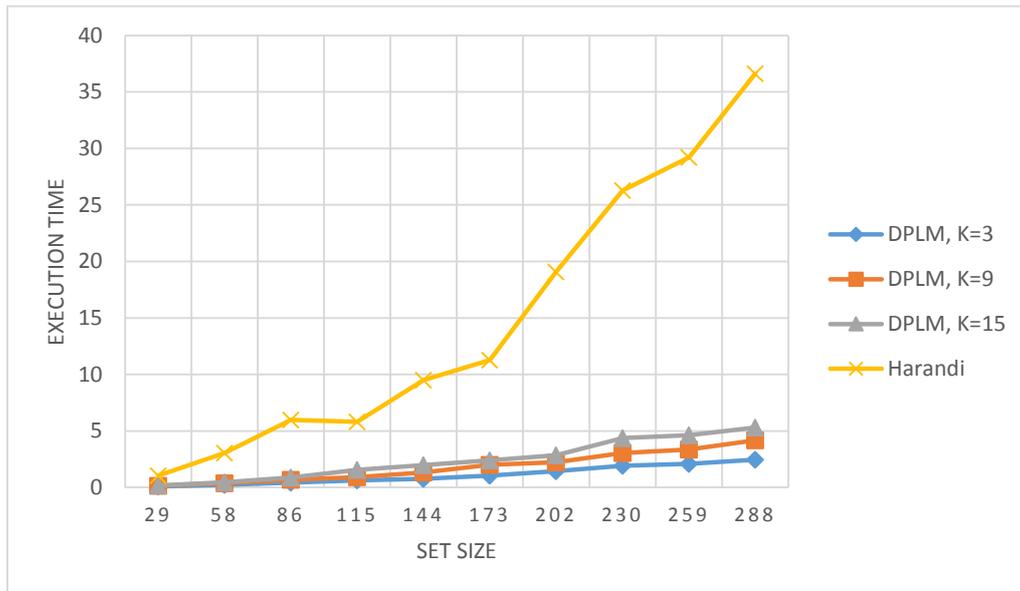

(a)

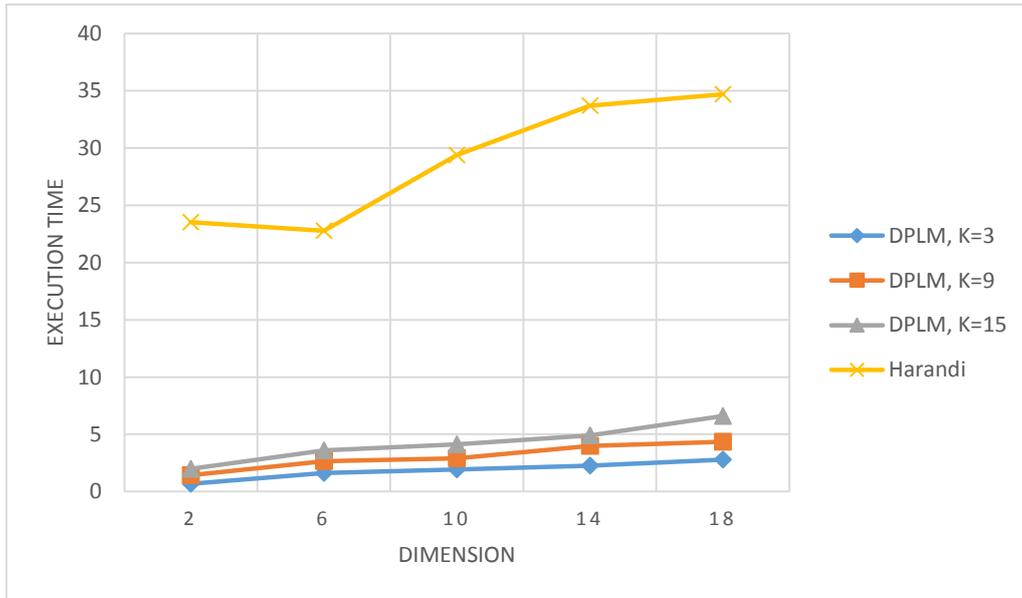

(b)

**Fig. 3.** Comparison of the execution time between DPLM and Harandi's method (a) execution time versus number of samples and (b) execution time versus number of dimensions.

## 5. Conclusion

We proposed a dimensionality reduction algorithm for the manifold of symmetric positive definite matrices, which preserves local structure of data by preserving distance to local means. This algorithm maps from a high dimensional SPD manifold to a lower dimensional one and can also use the label information in order to perform better in classification tasks. We compared our algorithm with the only similar method in the literature, i.e. Harandi's algorithm. We showed that our algorithm performs significantly better in terms of both execution time and accuracy. As our algorithm is linear in terms of number of samples, it can simply scale to larger number of samples without worrying about the execution time. We also proposed a time window and band-pass filter selection algorithm, which can find proper parameters for each subject and as we showed this would result in better performances in MI BCI application. We showed that using the proposed dimensionality reduction algorithm in conjunction with FGMDM as the classifier leads to the state of the art results on dataset IIa of BCI competitions IV.